\input amstex
\documentstyle{amsppt}
\topmatter
\title
A note on zeros of $L$-series of elliptic curves
\endtitle
\author
Xian-Jin Li
\endauthor
\abstract
In this note we study an analogy between a positive definite
quadratic form for elliptic curves over finite fields and a
positive definite quadratic form for elliptic curves over the
rational number field.  A question is posed of which an affirmative
answer would imply the analogue of the Riemann hypothesis for
elliptic curves over the rational number field.
\endabstract
\address
Department of Mathematics, Brigham Young University, Provo, Utah
84602 USA
\endaddress
\email
xianjin\@math.byu.edu
\endemail
\thanks
Research supported by National Security Agency H98230-06-1-0061
\endthanks
\endtopmatter
\document

\subheading{1.  Introduction}
  In the last entry of his diary, Gauss had noted the number
  $N$ of integer solutions of the congruence
  $$x^2y^2+x^2+y^2 \equiv 1 \mod p$$
  for a prime number $p\equiv 1 (\mod 4)$ is given by the formula
  $$N=p+1-(\pi+\bar \pi),$$
  where $\pi$ is defined by the decomposition
  $p=\pi\bar\pi$ in the ring $\Bbb Z[i]$.  It follows that
  $$|p+1-N|\leq 2\sqrt p. $$
  Gauss did not give a proof for his statement.
    In 1921 [7], G. Herglotz presented a proof of Gauss' statement
using complex multiplication.

   In 1933 H. Hasse [5] [6] used complex multiplication,
which is in essence the same method that G. Herglotz used,
to prove the analogue of the
Riemann hypothesis for all elliptic function fields.
 Let $\Bbb F_q$ be a finite field with $q$ elements,
and let  $E$ be an elliptic curve over $\Bbb F_q$.
 H. Hasse proved the inequality
$$|q+1-\#E(\Bbb F_q)|\leq 2\sqrt q, \tag 1.1$$
where $\#E(\Bbb F_q)$ is the number of $\Bbb F_q$-rational
points of the elliptic curve.

  In his book [4], H. Hasse said: ``The analogous procedure for
  the ordinary RH in algebraic number fields would be to
  get a best possible estimate of the prime number function
  $\Pi(x)$ or rather its analogue for prime divisors in algebraic
  number fields.  From Riemann's exact formula one can deduce that
  an estimate with error term $O(\sqrt x)$ or only
  $O(\sqrt x\log x)$ would establish $RH$.  Nobody however
  has been trying this most natural approach to the ordinary
  RH up to now!"

     Although the author does not understand the statement
  above,   he feels that it may indicate a method to attack the
  analogue of the   Riemann hypothesis for elliptic curves over the
  rational number field, which is the content of this note.

     We briefly review a proof of Hasse's inequality (1.1) in section 2.
  In section 3, an analogous procedure is presented for elliptic curves
  over the rational number field.  A question is posed of which an
affirmative answer would imply the analogue of the Riemann hypothesis
for elliptic curves over the rational number field.

     This note was initiated at the Workshop on
Zeta-Functions and Associated Riemann Hypotheses, New York University,
Manhattan, May 29 - June 1, 2002.  He wants to thank the American
Institute of Mathematics, Brian Conrey, and Peter Sarnak for the
financial support for him to attend the workshop.

\subheading{2.  L-series of elliptic curves over finite fields}
    Let $A$ be an abelian group.  A function $d: A\to\Bbb R$
is called a positive definite quadratic form if
$d(\alpha)=d(-\alpha)\geq 0$ for all $\alpha\in A$
with equality if and only if $\alpha=0$, and if the pairing
$$\langle\alpha,\beta\rangle=d(\alpha+\beta)-d(\alpha)-d(\beta)$$
is bilinear on $A$.
The degree map $\deg: \text{Hom}(E, E)\to\Bbb Z$ is a
positive definite quadratic form; see [10].  Denote
$$\langle\psi,\phi\rangle=\deg(\psi+\phi)-\deg(\psi)-\deg(\phi).$$
Since $\langle m\psi+n\phi, m\psi+n\phi\rangle\geq 0$
for all integers $m$ and $n$, we have the Cauchy-Schwartz inequality
$$|\deg(\psi-\phi)-\deg\psi-\deg\phi|\leq 2\sqrt{\deg \psi\deg\phi}
\tag 2.1$$
for all morphisms $\psi, \phi\in \text{Hom}(E,E)$.  Let
$$y^2+a_1xy+a_3y=x^3+a_2x^2+a_4x+a_6$$
be a Weierstrass equation for $E$ with coefficients in $\Bbb F_q$.  If
$$\phi: E\to E,\,\, (x,y)\to (x^q, y^q)$$
denotes the $q^{\text{th}}$-power Frobenius morphism of $E$,
then $\deg \phi=q$ and
$$\deg (1-\phi)=\#E(\Bbb F_q).$$
Thus, Hasse's inequality (1.1) follows from (2.1);
see [10].

     Let $a=q+1-\#E(\Bbb F_q)$.  The $L$-series of the elliptic
curve $E$ over $\Bbb F_q$ is given by
$$L_E(s)=1-aq^{-s}+q^{1-2s}.$$
 By Hasse's inequality we have
$$L_E(s)=(1-\alpha q^{-s})(1-\beta q^{-s})$$
with $\bar\alpha=\beta$ and $|\alpha|=|\beta|=\sqrt q$.
Hence, the analogue of the Riemann hypothesis for
L-series of elliptic curves over finite fields
follows from the Cauchy-Schwartz inequality (2.1).

\subheading{3.  L-series of elliptic curves over $\Bbb Q$}
   Let $E$ be an elliptic curve over the rational number field
$\Bbb Q$, and let $\hat h$ the canonical height on $E$. By
Artin's product formula, we have
$$\hat h(P)\geq 0$$
 for all points $P\in E(\bar{\Bbb Q})$, where $\bar{\Bbb Q}$
 is the algebraic closure of $\Bbb Q$ in the complex numbers.
 The N\'eron-Tate pairing on
$E$ is the bilinear form
$$\langle P, Q\rangle=\hat h(P+Q)-\hat h(P)-\hat h(Q)  $$
for all points $P, Q\in E(\bar{\Bbb Q})$; see [10].
Since $\langle m P+nQ, m P+nQ\rangle\geq 0$
for all integers $m$ and $n$, we have the Cauchy-Schwartz inequality
$$|\hat h(P-Q)-\hat h(P)-\hat h(Q)|
\leq 2\sqrt{\hat h(P)\hat h(Q)}\tag 3.1 $$
for all points $P,Q\in E(\bar{\Bbb Q})$.  Note that the fundamental
theorem of arithmetic is implicitly used in (3.1)
because of Artin's product formula.

   Let $N$ be the conductor of $E$.
For each prime $p$, we denote the reduction of $E$ at $p$ by
$\tilde E_p$.  Let
$$a_p=\cases p+1-\#\tilde E_p(\Bbb F_p), &\text{if $E$ has good reduction at $p$} \\
1, &\text{if $E$ has split multiplicative reduction at $p$} \\
-1, &\text{if $E$ has non-split multiplicative reduction at $p$} \\
0, &\text{if $E$ has additive reduction at $p$.}\endcases$$
The $L$-series associated to the elliptic curve
$E$ is defined by the Euler product:
$$L_E(s)=\prod_{p\nmid N}(1-a_pp^{-s}+p^{1-2s})^{-1}
\prod_{p|N}(1-a_p p^{-s})^{-1} $$
for $\Re s>3/2$.

   By results in [3] [9] and the Shimura-Taniyama conjecture that
is a theorem proved in [2] [11],
there is a normalized Hecke eigenform $f$ of weight $2$ and level $N$
and of rational Fourier coefficients such that $L_f=L_E$.
  Hence $L_E(s)$ has an analytic
continuation to the complex plane and satisfies the
functional equation
$$\Gamma(s)L_E(s)
=wN^{1-s}(2\pi)^{2s-2}\Gamma(2-s)L_E(2-s), $$
where $w=(-1)^r$ with $r$ being the vanishing order of $\xi_f(s)$ at
$s=1/2$; see [1].

 Let
$$\xi_E(s)=N^{s/2}(2\pi)^{-s}\Gamma({1\over 2}+s)L_E({1\over 2}+s).$$
Then $\xi_E(s)$ is an entire function and satisfies the functional identity
$$\xi_E(s)=w\xi_E(1-s).$$
Let $\varphi(z)=\xi_E(1/(1-z))$, and let
$${\varphi^\prime(z)\over\varphi(z)}
=\sum_{n=0}^\infty\lambda_E(n+1)z^n.\tag 3.2$$
Similarly as in [8], we can show that all zeros of $\xi_E(s)$
lie on the line $\Re s=1/2$ if and only if $\lambda_E(n)\geq 0$
for $n=1,2,\cdots$.

 On the other hand, to show that all zeros of $\xi_E(s)$
 lie on the line $\Re s=1/2$
it suffices to obtain upper bounds for the $\lambda_E(n)$'s
that would give the convergence of the series (3.2) for
$|z|<1$.  Studying H. Hasse's proof of the Riemann
hypothesis for elliptic curves over finite fields leads
 the author to believe that the following question should have an
affirmative answer.

\proclaim{Question}  Is it possible that there exist points
$P, Q\in E$ such that certain upper bounds for the
$\lambda_E(n)$'s, which would give the convergence of the series
(3.2) in the unit disk $|z|<1$, can be obtained by using
the Cauchy-Schwartz inequality (3.1)?
\endproclaim


\Refs
\ref
  \no 1
  \by A. O. L. Atkin and J. Lehner
  \paper Hecke operators on $\Gamma \sb{0}(m)$
  \jour Math. Ann.
  \vol 185
  \yr 1970
  \pages 134--160
\endref
\ref
  \no 2
  \by C. Breuil, B. Conrad, F. Diamond, and R. Taylor
  \paper On the modularity of elliptic curves over $\Bbb Q$:
   wild 3-adic exercises
  \jour J. Amer. Math. Soc.
  \vol 14
  \yr 2001
  \pages 843--939
\endref
\ref
  \no 3
  \by H. Carayol
  \paper Sur les repr\'esentations $\ell$-adiques associ\'ees aux formes
modulaires de Hilbert
  \jour Ann. Sci. \'Ecole Norm. Sup.
  \vol 19
  \yr 1986
  \pages 409--468
\endref
\ref
  \no 4
  \by H. Hasse
  \book The Riemann Hypothesis in Algebraic Function Fields over
        a finite constants field
  \publ Pennsylvania State University
  \publaddr University Park
  \yr 1968
\endref
\ref
  \no 5
  \by H. Hasse
  \paper Beweis des  Analogons der Riemannschen Vermutung f\"ur
  die Artinschen und F. K. Schmidtschen Kongruenzzetafunktionen
  in gewissen elliptischen F\"allen. Vorl\"aufige Mitteilung
  \publ in ``Helmut Hasse Mathematische Abhandlungen,"
   Band 2, deGruyter, 1975
   \pages 85--94
\endref
\ref
  \no 6
  \by H. Hasse
  \paper Modular functions and elliptic curves over finite fields
  \publ in ``Helmut Hasse Mathematische Abhandlungen,"
   Band 2, deGruyter, 1975
   \pages 351--369
\endref
\ref
  \no 7
  \by G. Herglotz
  \paper  Zur letzten Eintragung im Gausschen Tagebuch
  \jour Ber. Verhandl. S\"achs. Akad. Wiss. Math.-Phys. Kl.
  \vol 73
  \yr 1921
  \pages 271--276
\endref
\ref
  \no 8
  \by Xian-Jin Li
  \paper The positivity of a sequence of numbers and the Riemann hypothesis
  \jour J. Number Theory
  \vol 65
  \yr 1997
  \pages 325--333
\endref
\ref
  \no 9
  \by G. Shimura
  \book Introduction to the Arithmetic Theory of Automorphic Functions
  \publ Princeton University Press
  \yr 1971
\endref
\ref
  \no 10
  \by J. H. Silverman
  \book  The Arithmetic of Elliptic Curves
  \publ Springer-Verlag
  \publaddr New York
  \yr 1986
\endref
\ref
  \no 11
  \by A. Wiles
  \paper Modular elliptic curves and Fermat's last theorem
  \jour  Ann. of Math.
  \vol 141
  \yr 1995
  \pages 443--551
\endref
\endRefs
\enddocument